\documentclass[11pt,reqno]{amsart}

\usepackage[T1]{fontenc}
\usepackage{mathptmx}
\usepackage{microtype}
\usepackage{amsmath,amssymb,amsthm,mathtools,mathrsfs}
\usepackage{graphicx}
\usepackage{comment}
\usepackage{enumitem}
\usepackage[colorlinks=true,linkcolor=blue,citecolor=blue,urlcolor=blue]{hyperref}
\usepackage[nameinlink,capitalise]{cleveref}
\usepackage{geometry}
\mathtoolsset{showonlyrefs}

\newcommand{\C}{\mathbb{C}}

\newcommand{\Z}{\mathbb{Z}}
\newcommand{\R}{\mathbb{R}}

\newcommand{\calA}{\mathcal{A}}

\newcommand{\calT}{\mathcal{T}}

\newcommand{\e}{\mathrm{e}}

\numberwithin{equation}{section}
\newtheorem{thm}{Theorem}[section]

\newtheorem{prop}[thm]{Proposition}

\theoremstyle{definition}

\theoremstyle{remark}
\newtheorem{rem}[thm]{Remark}

\title[McShane's identity and the Brownian loop measure]
{McShane's identity and the Brownian loop measure}
\begin{document}
\author[Hedenmalm]{Haakan Hedenmalm}

\address{Department of Mathematics
and Computer Science, St Petersburg State University, St Petersburg, Russia.
\\
Beijing Institute for Mathematical Sciences and Applications,
Huairou District, Beijing 101408, China
\\
Department of Mathematics and Statistics, University of Reading, U. K.
\\
Department of Mathematics, University of Aveiro, Portugal.
}
\email{haakan00@gmail.com}



\date{\today}

\begin{abstract}
We combine McShane's identity for a once-punctured torus with the
length-spectrum identity of Wang and Xue. An exact identity in each
primitive homology class yields a finite limit for the difference of
two divergent length-spectrum sums over self-intersecting hyperbolic
geodesics. The limit is $\frac18$ and is equivalent to McShane's identity.
We also derive higher-iteration identities governed by integer divisors.
\end{abstract}

\maketitle

Wang and Xue explicitly ask, immediately after Remark~1.4 of
their paper~\cite{WangXue}, whether their length-spectrum identity
is related to McShane's identity for the once-punctured torus.
The present note gives an explicit answer to this question by
comparing the Wang--Xue formula for a primitive flat homology class
with the formula for its double.

The purpose of this note is to express McShane's identity as a
convergent difference of length-spectrum sums, taken over the same
finite collections of primitive homology classes. While the two
individual sums diverge, their difference converges to $\frac18$,
independently of the exhaustion. The distinction between geometric
primitivity and primitive homology is essential. The last section
compares higher multiples via the Wang--Xue identity.

%\tableofcontents

\section{McShane's identity}
We let $\calT$ be a torus with complex structure, which we may model as
$\C/\Lambda$, where $\Lambda$ is a discrete translation group with two
$\R$-linearly independent generators. We take a point $\zeta_0\in\calT$ and
consider the \emph{punctured torus} $\calT':=\calT\setminus\{\zeta_0\}$.
The torus $\calT$ is given a flat metric (locally, it is the metric of
$\C$ carried over to $\C/\Lambda$). On the other hand, the punctured torus
$\calT'$ is given the unique complete hyperbolic metric (with constant
negative curvature $-1$) which is conformal to the complex structure. We let
$\ell_\gamma(\calT)$ denote the length of the curve $\gamma$ with respect
to the flat metric of $\calT$, and, likewise, $\ell_{\gamma'}(\calT')$ is
the length of the curve $\gamma'$ with respect to the hyperbolic metric
of $\calT'$. 

\begin{thm}
{\rm (McShane's identity)} We have that
\[
\sum_{\gamma'} \frac{1}{\e^{\ell_{\gamma'}(\calT')}+1}=\frac12, 
\]
where the sum runs over all (unoriented) simple closed geodesics $\gamma'$
on $\calT'$. 
\end{thm}

The relevant reference is \cite{McShane}.

\begin{rem}
The expression
\[
\frac{1}{\e^{\ell_{\gamma'}(\calT')}+1}
\]
bears a strong resemblance with Fermi-Dirac statistics (see, e.g.,
\cite{LandauLifshits}).
After all, in Fermi-Dirac statistics the expected number of fermions
in state $i$ is given by
\[
\bar n_i=\frac{1}{\e^{(\epsilon_i-\mu)/(k_B T)}+1}
\]
Here, $\epsilon_i$ is the energy of the state $i$, $\mu$ is the chemical
potential, $k_B$ is Boltzmann's constant, and $T$ the absolute temperature.
Suppose that our particles in the Fermi-Dirac model consist of homotopy classes
of loops on $\calT'$, and that the energy $\epsilon$ of a given homotopy class
is such that
\[
\frac{\epsilon-\mu}{k_B T}=\ell_{\gamma'}(\calT'),
\]
where $\gamma'$ is the minimal length geodesic in the given homotopy class,
then McShane's identity becomes an assertion about the expected number of
homotopy classes with $\epsilon>\mu$ in the Fermi-Dirac model that are 
homotopically equivalent to a simple geodesic loop on $\calT'$.
\end{rem}

\section{The work of Wang-Xue based on the Brownian
loop measure}

In the recent work \cite{WangXue}, Yilin Wang and Yuhao Xue study the
Brownian loop measure introduced by Lawler and Werner \cite{LawlerWerner} and establish
a remarkable identity of surface length spectra on a Riemann surface
$\mathcal{X}$
and on $\mathcal{X}'$, which is $\mathcal{X}$ but with finitely many points
removed. Applied in our context of a flat torus $\calT$ and hyperbolic torus
$\calT'$ with a single point removed, they obtain the following. 

\begin{thm}
{\rm (Wang-Xue)} Let $\gamma$ be a closed geodesic on the flat torus $\calT$.
We then have that
\[  
\frac{\mathrm{Area}(\calT)}{\pi\ell_\gamma(\calT)^2}=\sum_{\gamma'\simeq_\calT \gamma}
\frac{1}{m(\gamma')}\,\frac{1}{\e^{\ell_{\gamma'}(\calT')}-1},
\]  
where $\gamma'$ runs over the closed geodesics on $\calT'$ that are homotopic
to $\gamma$. The number $m(\gamma')\in\{1,2,3,\ldots\}$ is the iteration
number of the closed geodesic $\gamma'$, which is a divisor of the iteration
number $m(\gamma)$ of the geodesic of the flat torus $\calT$. 
\end{thm}  

This is Theorem 4.3 of \cite{WangXue}. In this theorem the free
homotopy classes are taken with orientation. We apply it first to each
oriented class and then identify a class with its reverse. Both sides are
unchanged by reversal, including the iteration number, so this convention
does not introduce a factor of two. All subsequent sums over geodesics
are unoriented unless explicitly stated otherwise.

\begin{rem}
%The expression
%\[
%\frac{1}{\e^{\ell_{\gamma'}(\calT')}-1}
%\]
%bears a strong resemblance with Bose-Einstein statistics.
We write the identity of the Wang-Xue theorem in the form
\[  
\frac{m(\gamma)\,\mathrm{Area}(\calT)}{\pi\ell_\gamma(\calT)^2}
=\sum_{\gamma'\simeq_\calT \gamma}
\frac{m(\gamma)/m(\gamma')}{\e^{\ell_{\gamma'}(\calT')}-1},
\] 
and observe that the ratios $m(\gamma)/m(\gamma')$ are positive integers.
In contrast with the apparent Fermi-Dirac statistics conjured up by the
McShane theorem, here the expressions suggest
instead Bose-Einstein statistics. After all, in Bose-Einstein statistics
(see, e.g., \cite{LandauLifshits}), the expected number of bosons in state
$i$ is given by
\[
\bar n_i=\frac{d_i}{\e^{(\epsilon_i-\mu)/(k_B T)}-1},
\]
where $d_i$ is the degeneracy of energy level $i$, $\epsilon_i$
is its energy, $\mu$ is the chemical potential, $k_B$ is Boltzmann's
constant, and $T$ is the absolute temperature. The correspondence
$($energy minus chemical potential$)/(k_BT)=\ell_{\gamma'}(\calT')$
identifies the kernels exactly. The integer coefficient
$m(\gamma)/m(\gamma')$ is supplied by the iteration structure of
Wang--Xue's identity and may be viewed as a \emph{formal degeneracy
factor}. It is not asserted to count independent quantum states.
\end{rem}

\section{Bosonic interpretation of McShane's identity}

Next, let $\gamma_0$ be a simple closed geodesic on the flat torus $\calT$,
which means that it is closed, has finite length, and no self-intersection.
It can be modelled as a straight line segment in $\C/\Lambda$, such that
the starting point is equivalent to the endpoint modulo the lattice group
$\Lambda$, and that no other point inbetween is equivalent to the starting
point. We also need $\gamma_0^2$, which is just the loop $\gamma_0$ traversed
twice.
We first apply the Wang-Xue theorem with $\gamma=\gamma_0$:
\begin{equation}
  \frac{\mathrm{Area}(\calT)}{\pi\ell_{\gamma_0}(\calT)^2}=
  \sum_{\gamma'\simeq_\calT \gamma_0}
\frac{1}{\e^{\ell_{\gamma'}(\calT')}-1},
\label{eq:WX1}
\end{equation}
since the iteration number $m(\gamma')=1$ for any such loop $\gamma'$.
Next, we apply the Wang-Xue theorem with $\gamma=\gamma_0^2$:
\begin{equation}
\frac{\mathrm{Area}(\calT)}{\pi\ell_{\gamma_0^2}(\calT)^2}
=\sum_{\gamma''\simeq_\calT \gamma_0^2}
\frac{1}{m(\gamma'')}\,\frac{1}{\e^{\ell_{\gamma''}(\calT')}-1},
\label{eq:WX2}
\end{equation}
where we realize that the iteration number of $\gamma''$ is either
$m(\gamma'')=1$ or $m(\gamma'')=2$. If $m(\gamma'')=2$, then
$\gamma''=(\gamma')^2$, where $\gamma'$ is a primitive closed
geodesic on $\calT'$. It need not be simple. The geodesic $\gamma'$
is homotopic to $\gamma_0$ on the flat torus $\calT$, up to orientation.
The length of such a geodesic is $\ell_{\gamma''}(\calT')=
2\ell_{\gamma'}(\calT')$, and likewise on the flat torus, we have 
$\ell_{\gamma_0^2}(\calT)=2\ell_{\gamma_0}(\calT)$. It now follows that we
may express \eqref{eq:WX2} in the form
\begin{equation}
\frac{\mathrm{Area}(\calT)}{4\pi\ell_{\gamma_0}(\calT)^2}
=\sum_{\gamma''\simeq_\calT \gamma_0^2
,\,m(\gamma'')=1}
\frac{1}{\e^{\ell_{\gamma''}(\calT')}-1}
+\sum_{\gamma'\simeq_\calT \gamma_0}
\frac12\,\frac{1}{\e^{2\ell_{\gamma'}(\calT')}-1}.
\label{eq:WX2.1}
\end{equation}
Multiplying both sides by $4$ yields
\begin{equation}
\frac{\mathrm{Area}(\calT)}{\pi\ell_{\gamma_0}(\calT)^2}
=\sum_{\gamma''\simeq_\calT \gamma_0^2
,\,m(\gamma'')=1}
\frac{4}{\e^{\ell_{\gamma''}(\calT')}-1}
+\sum_{\gamma'\simeq_\calT \gamma_0}
\frac{2}{\e^{2\ell_{\gamma'}(\calT')}-1}.
\label{eq:WX2.2}
\end{equation}
The left-hand sides of \eqref{eq:WX1} and \eqref{eq:WX2.2} are the same, so
their right-hand sides must be equal too:
\begin{equation}
\sum_{\gamma'\simeq_\calT \gamma_0}
\frac{1}{\e^{\ell_{\gamma'}(\calT')}-1}=\sum_{\gamma''\simeq_\calT \gamma_0^2
,\,m(\gamma'')=1}
\frac{4}{\e^{\ell_{\gamma''}(\calT')}-1}
+\sum_{\gamma'\simeq_\calT \gamma_0}
\frac{2}{\e^{2\ell_{\gamma'}(\calT')}-1}.
\label{eq:WX3}
\end{equation}
Reshuffling terms in \eqref{eq:WX3}, we find that
\begin{multline}
\sum_{\gamma''\simeq_\calT \gamma_0^2
,\,m(\gamma'')=1}
\frac{4}{\e^{\ell_{\gamma''}(\calT')}-1}
= \sum_{\gamma'\simeq_\calT \gamma_0}
\frac{1}{\e^{\ell_{\gamma'}(\calT')}-1}-\sum_{\gamma'\simeq_\calT \gamma_0}
\frac{2}{\e^{2\ell_{\gamma'}(\calT')}-1}
\\
=\sum_{\gamma'\simeq_\calT \gamma_0}\bigg(
\frac{1}{\e^{\ell_{\gamma'}(\calT')}-1}-\frac{2}{\e^{2\ell_{\gamma'}(\calT')}-1}\bigg)
=\sum_{\gamma'\simeq_\calT \gamma_0}\frac{1}{\e^{\ell_{\gamma'}(\calT')}+1}.
\label{eq:WX3.1}
\end{multline}
We now sum the classwise identity over a finite collection of primitive
homology classes. A nonzero class in $H_1(\calT,\Z)$ is called
\emph{primitive} if it is not an integer multiple by an integer of
absolute value greater than $1$ of another integral class. A closed
hyperbolic geodesic is called \emph{primitive} if it is not a proper
iterate. These are distinct notions. We identify classes differing by
sign since the geodesics are unoriented.

Let $\calA$ be the set of primitive nonzero homology classes of
$\calT$, modulo sign. For a finite subset $F\subset\calA$, put
\begin{align*}
 B(F)&=\sum_{a\in F}\ \sum_{\substack{\gamma''\,\text{\rm primitive}\\
 [\gamma'']=2a}}\frac{1}{\e^{\ell_{\gamma''}(\calT')}-1},\\
 E(F)&=\sum_{a\in F}\ \sum_{\substack{\gamma'\,\text{\rm self-intersecting}\\
 [\gamma']=a}}\frac{1}{\e^{\ell_{\gamma'}(\calT')}+1}.
\end{align*}
The sums in each fixed homology class converge by the Wang--Xue
identity. Every simple closed hyperbolic geodesic on the once-punctured
torus represents a primitive nonzero homology class on the filled torus.
Indeed, it is nonperipheral and nonseparating. Conversely, a closed
geodesic with primitive homology is primitive as a closed curve, though
it need not be simple.

\begin{thm}\label{thm:bosonic}
For every finite subset $F\subset\calA$, we have the exact identity
\begin{equation}\label{eq:finite-bosonic}
 B(F)-\frac14E(F)
 =\frac14\sum_{\substack{\gamma\,\text{\rm simple}\\{}[\gamma]\in F}}
 \frac{1}{\e^{\ell_\gamma(\calT')}+1}.
\end{equation}
Consequently, for any increasing sequence of finite sets
$F_j\subset\calA$ with $\bigcup_jF_j=\calA$,
\begin{equation}\label{eq:bosonicMcShane}
 \lim_{j\to\infty}\left(B(F_j)-\frac14E(F_j)\right)=\frac18.
\end{equation}
The limit statement is equivalent to McShane's identity, in the
presence of the Wang--Xue identity.
\end{thm}

\begin{proof}
For every $a\in\calA$, the comparison of the Wang--Xue identities
for $a$ and $2a$ in \eqref{eq:WX3.1} gives
\[
 4\sum_{\substack{\gamma''\,\text{\rm primitive}\\{}[\gamma'']=2a}}
 \frac1{\e^{\ell_{\gamma''}(\calT')}-1}
 =\sum_{[\gamma']=a}\frac1{\e^{\ell_{\gamma'}(\calT')}+1}.
\]
Each geodesic in the sum on the right has primitive homology, hence is
primitive as a closed geodesic. Partitioning this sum into the simple
and self-intersecting geodesics, and summing over $a\in F$, proves
\eqref{eq:finite-bosonic}. McShane's identity says that the sum over
all simple unoriented closed geodesics of the kernel
$(\e^{\ell}+1)^{-1}$ equals $\frac12$. Its summands are positive,
so taking the increasing limit yields \eqref{eq:bosonicMcShane}.
Conversely, \eqref{eq:finite-bosonic} holds independently of McShane's
identity. If \eqref{eq:bosonicMcShane} holds for such an exhaustion,
then passage to the limit in \eqref{eq:finite-bosonic} yields
McShane's identity.
\end{proof}

\begin{rem}
The individual terms $B(F_j)$ and $E(F_j)$ diverge as $j\to\infty$.
Indeed, Wang--Xue gives
\[
 \sum_{[\gamma]=a}\frac1{\e^{\ell_\gamma}-1}
 =\frac{\mathrm{Area}(\calT)}{\pi L_a^2},
\]
where $L_a$ is the length of the primitive flat representative.
The sum of $L_a^{-2}$ over primitive nonzero lattice classes diverges.
As the lengths of closed hyperbolic geodesics are uniformly bounded
away from zero, the corresponding fermionic and bosonic kernels are
uniformly comparable. The cancellation in \eqref{eq:finite-bosonic}
is therefore essential. No subtraction of untruncated infinite series
is intended. A natural choice of exhaustion is
$F_R=\{a\in\calA:L_a\le R\}$.
\end{rem}

The logarithmic divergence of these two sums can be quantified. This
also shows directly how their leading terms cancel.

\begin{prop}\label{prop:logarithmic}
For $F_R=\{a\in\calA:L_a\le R\}$, put $B_R=B(F_R)$ and
$E_R=E(F_R)$. As $R\to+\infty$,
\begin{align*}
 B_R&=\frac{3}{2\pi^2}\log R+\mathrm{O}(1),\\
 E_R&=\frac{6}{\pi^2}\log R+\mathrm{O}(1).
\end{align*}
In particular, the logarithmic terms cancel in
$B_R-\frac14 E_R$, whose limit is $\frac18$.
\end{prop}

\begin{proof}
Write $A=\mathrm{Area}(\calT)$. Primitive lattice-point counting,
with the identification of $a$ and $-a$, gives
\[
\sum_{a\in F_R}\frac{A}{\pi L_a^2}
=\frac6{\pi^2}\log R+\mathrm{O}(1).
\]
Indeed, the number of primitive nonzero vectors in the flat lattice,
modulo sign, with length at most $R$ is
$3R^2/(\pi A)+\mathrm{O}(R\log(2+R))$.
Partial summation yields the displayed logarithmic asymptotic.

For $x$ bounded away from zero, the difference between
$(\e^x-1)^{-1}$ and $(\e^x+1)^{-1}$ is
$2/(\e^{2x}-1)=\mathrm{O}(\e^{-2x})$.
The sum of $\e^{-2\ell_\gamma(\calT')}$ over primitive closed
hyperbolic geodesics converges, by the standard bound
$\#\{\gamma: \ell_\gamma(\calT')\le t\}=\mathrm{O}(\e^t)$
for a fixed finite-area hyperbolic surface. Hence
\[
\sum_{a\in F_R}\sum_{[\gamma]=a}
\frac1{\e^{\ell_\gamma(\calT')}+1}
=\frac6{\pi^2}\log R+\mathrm{O}(1).
\]
The sum over simple closed geodesics of this kernel remains bounded
by McShane's identity. Removing those terms gives the asserted
asymptotic for $E_R$. Finally the classwise identity in the proof of
Theorem~\ref{thm:bosonic} gives
$4B_R=E_R+\mathrm{O}(1)$, proving the first asymptotic as well.
\end{proof}

\section{Higher iterations and arithmetic identities}

The preceding comparison of a primitive flat geodesic with its double
extends to all positive iteration numbers. Fix a primitive nonzero class
$a\in H_1(\calT,\Z)$, understood up to sign, and let $L_a$ be the length
of the corresponding simple flat geodesic. Put
\[
C_a=\frac{\mathrm{Area}(\calT)}{\pi L_a^2},\qquad
S_k(d;a)=\sum_{\substack{\gamma'\,\text{\rm primitive}\\
{}[\gamma']=ka}}
\frac{1}{\e^{d\ell_{\gamma'}(\calT')}-1},
\qquad k,d\ge1.
\]
Here and below, $[\gamma']=ka$ is interpreted up to sign. The primitivity
condition on $\gamma'$ refers to geometric iteration, not to homology.

\begin{prop}
For every positive integer $n$, we have
\begin{equation}
\frac{C_a}{n^2}=\sum_{d\mid n}\frac1d S_{n/d}(d;a).
\label{eq:divisoridentity}
\end{equation}
\end{prop}

\begin{proof}
Apply the Wang-Xue identity to the $n$-fold iterate of the simple flat
geodesic representing $a$. A hyperbolic geodesic in its filled-torus
homotopy class can be written uniquely as the $d$-fold iterate of a
primitive hyperbolic geodesic $\gamma'$. Its homology is then
$d[\gamma']=na$. Since $a$ is primitive, this implies $d\mid n$ and
$[\gamma']=(n/d)a$, up to sign. The iteration factor in the Wang-Xue
formula is $1/d$, while the length is $d\ell_{\gamma'}(\calT')$.
Since the flat length is $nL_a$, the asserted formula follows.
\end{proof}

For $n=1$, formula \eqref{eq:divisoridentity} gives $S_1(1;a)=C_a$.
For a prime $p$, its only divisors are $1$ and $p$. We obtain the
following extension of \eqref{eq:WX3.1}.

\begin{thm}
For every prime $p$ and primitive nonzero class $a\in H_1(\calT,\Z)$,
\begin{equation}
p^2\sum_{\substack{\gamma'\,\text{\rm primitive}\\{}[\gamma']=pa}}
\frac{1}{\e^{\ell_{\gamma'}(\calT')}-1}
=\sum_{[\gamma']=a} K_p\big(\ell_{\gamma'}(\calT')\big),
\label{eq:primeidentity}
\end{equation}
where the geodesics in the sum on the right are necessarily primitive and
\begin{equation}
K_p(x)=\frac1{\e^x-1}-\frac{p}{\e^{px}-1},\qquad x>0.
\label{eq:primekernel}
\end{equation}
The kernel $K_p$ is strictly positive. Moreover,
\[
\lim_{x\to0^+}K_p(x)=\frac{p-1}{2},\qquad
K_p(x)\sim\e^{-x}\quad(x\to+\infty).
\]
\end{thm}

\begin{proof}
The divisor identity for $n=p$ reads
$C_a/p^2=S_p(1;a)+p^{-1}S_1(p;a)$.
Subtracting $pS_1(p;a)$ from $S_1(1;a)=C_a$ gives
\eqref{eq:primeidentity}. For $t=\e^x>1$, the inequality
$1+t+\cdots+t^{p-1}>p$ shows that $K_p(x)>0$.
The limit at zero and the asymptotic at infinity follow directly from
\eqref{eq:primekernel}.
\end{proof}

For $p=2$ we recover $K_2(x)=(\e^x+1)^{-1}$, the fermionic kernel in
McShane's identity. For $p=3$, the new kernel takes the form
\[
K_3(x)=\frac{\e^x+2}{\e^{2x}+\e^x+1}.
\]
Thus the case $p=2$ is the first member of a family of positive
length-spectrum identities. The special constant in Theorem 3.2 comes
from applying McShane's identity specifically to the kernel $K_2$.

\begin{rem}[Gentile statistics]
The kernel in \eqref{eq:primekernel} has an exact interpretation in
finite-occupancy quantum statistics. For every integer $p\ge2$, put
\[
Z_p(x)=\sum_{j=0}^{p-1}\e^{-jx}
       =\frac{1-\e^{-px}}{1-\e^{-x}},\qquad x>0.
\]
Then
\begin{equation}
K_p(x)=-\frac{\mathrm{d}}{\mathrm{d}x}\log Z_p(x)
=\frac{\sum_{j=0}^{p-1}j\e^{-jx}}
       {\sum_{j=0}^{p-1}\e^{-jx}}.
\label{eq:gentile}
\end{equation}
Thus $K_p(x)$ is the mean occupation number of one level when the
allowed occupancies are $0,1,\ldots,p-1$, with $x$ interpreted as the
dimensionless energy relative to the chemical potential. This is the
Gentile finite-occupancy model. At $p=2$ it is the Fermi--Dirac
occupation number, whereas $Z_p(x)$ and $K_p(x)$ tend as $p\to\infty$
to the Bose--Einstein partition function and occupation number for
fixed $x>0$. The interpretation is an exact identity of kernels. It
does not identify the hyperbolic geodesics themselves with physical
particles.
\end{rem}

For composite values of $n$, intermediate iteration numbers occur.
For instance, $n=4$ gives
\begin{equation}
16S_4(1;a)=S_1(1;a)-8S_2(2;a)-4S_1(4;a).
\label{eq:fourfold}
\end{equation}
This illustrates the divisor arithmetic that is absent from the prime
case.

\begin{rem}
The prime kernels admit a roots-of-unity decomposition. If
$\omega=\e^{2\pi\mathrm{i}/p}$, partial fractions give
\[
\frac{p}{t^p-1}=\sum_{j=0}^{p-1}\frac{\omega^j}{t-\omega^j},
\qquad
K_p(x)=-\sum_{j=1}^{p-1}\frac{\omega^j}{\e^x-\omega^j}.
\]
For $p=2$ this reduces to the single real term $(\e^x+1)^{-1}$.
For larger primes, the terms combine in complex-conjugate pairs to
produce the real, strictly positive kernel $K_p$.
\end{rem}

\begin{comment}
\subsection*{Acknowledgements}
The author acknowledges support from Vetenskapsrådet (VR grant 2020-03733),
the Leverhulme trust (grant VP1-2020-007), by grant 075-15-2021-602
of the Government of the Russian Federation for the state support of
scientific research, carried out under the supervision of leading
scientists, and by the St. Petersburg Leonhard Euler International
Mathematical Institute, under grant 075-15-2022-287. 
\end{comment}

\end{document}